\documentclass[12pt]{amsart}
\usepackage{amssymb,amsmath,amsthm}
\usepackage{hyperref}

\makeatletter\@addtoreset{equation}{section} \makeatother

\newcommand{\nad}[2]{\genfrac{}{}{0pt}{}{#1}{#2}}
\newcommand{\p}{\partial}

\renewcommand{\a}{\alpha}

\def \g {\gamma}
\def \a {\alpha}
\def \b {\beta}

\def \C {\mathbb{C}}

\def \Z {\mathbb{Z}}
\def \N {\mathbb{N}}

\def \Sympl {\mathfrak{S}}

\newtheorem{theorem}[equation]{Theorem}
\newtheorem{proposition}[equation]{Proposition}
\newtheorem{lemma}[equation]{Lemma}
\newtheorem{corollary}[equation]{Corollary}

\theoremstyle{remark}
\newtheorem{remark}[equation]{Remark}

\theoremstyle{definition}
\newtheorem{example}[equation]{Example}
\newtheorem{definition}[equation]{Definition}

\newcommand{\Diag}{\mathop{\mathrm{Diag}}\nolimits}

\title{Singular polynomials from orbit spaces}

\author{Misha\,Feigin, Alexey\,Silantyev}

\newcommand{\wt}{\widetilde}
\newcommand{\wh}{\widehat}

\newcommand{\Cstar}{{\mathop{C}\limits^*}\vphantom{|}}

\begin{document}

\begin{abstract}We consider the polynomial representation  $S(V^*)$ of the rational Cherednik algebra $H_c(W)$ associated to a finite Coxeter group $W$ at constant parameter
$c$.
We show that for any degree $d$ of $W$ and $m\in \N$ the space $S(V^*)$ contains a single copy of the reflection representation $V$ of $W$ spanned by the homogeneous singular polynomials of degree $d-1+h m$, where $h$ is the Coxeter number of $W$; these polynomials generate an $H_c(W)$ submodule with the parameter $c=(d-1)/h+m$. We express these singular polynomials through the Saito polynomials that are flat coordinates of the Saito metric on the orbit space $V/W$. We also show that this exhausts all the singular polynomials in the isotypic component of the reflection representation $V$ for any constant parameter $c$.

\vspace{6mm}

 \end{abstract}

\address{School of Mathematics and Statistics, University of Glasgow, 15 University Gardens,
Glasgow G12 8QW, UK}

\address{Graduate School of Mathematical Sciences, University of Tokyo, Komaba, Tokyo 153-8914, Japan}

\email{misha.feigin@glasgow.ac.uk, aleksejsilantjev@gmail.com}

\maketitle

\section{Introduction}

In this paper we relate two remarkable constructions associated with a finite Coxeter group $W$. The first one is the Frobenius manifold structure on the space of orbits of $W$ acting in its reflection representation $V$~\cite{Dubrovin_lect}. The key ingredient here is the Saito flat metric on the orbit space $V/W$~\cite{SYS}. This metric is defined as a Lie derivative of the standard contravariant (Arnold) metric. The flat coordinates form a distinguished basis in the ring of invariant polynomials $S(V^*)^W$. This basis is now known explicitly for all irreducible groups $W$. All the cases except $W$ of type $E_7$, $E_8$ were covered in the original paper~\cite{SYS}. The flat coordinates in the latter two cases were found recently both in~\cite{Abriani} and in~\cite{Talamini}.

The other famous construction associated with the group $W$ is the rational Cherednik algebra $H_c(W)$~\cite{EG}. It depends on the $W$-invariant function $c$ on the set of reflections of $W$ which we assume to be constant. The key ingredient here is the Dunkl operator $\nabla_\zeta$, $\zeta \in V$, which acts in the ring of polynomials as a differential-reflection operator~\cite{Dunkl}. For particular values of $c$ the polynomial representation $S(V^*)$ has non-trivial submodules $M$.  These values were completely determined by Dunkl, de Jeu and Opdam  in~\cite{DJO} where it was shown that non-trivial submodules exist if and only if $c$ is a non-integer number of the form $c=l/d$ where $d$ is one of the degrees of the Coxeter group $W$ and $l \in \Z_{>0}$. The lowest homogeneous component $M_0$
of $M$ consists of so-called singular polynomials~\cite{DJO} which are annihilated by Dunkl operators $\nabla_\zeta$ for any $\zeta \in V$. All singular polynomials were found by Dunkl when $W$ has type $A$~\cite{DunklAsing1, DunklAsing2}. Further, it was established in~\cite{ESG} that in this case any submodule $M$ is generated by its lowest homogeneous component $M_0$. In general the structure of submodules of $S(V^*)$ and the corresponding singular polynomials are not known.  Some singular polynomials for the classical groups $W$ and for the icosahedral group were determined in~\cite{EtChm} and~\cite{Du2} (see also~\cite{Balag}) respectively, while dihedral case was fully studied in~\cite{DJO} (see also~\cite{Chm}).

In the paper we study singular polynomials that belong to the isotypic component of the reflection representation $V$ of the Coxeter group $W$. The existence of such singular polynomials is known for the Weyl groups  when $c=r/h$ where $h$ is the Coxeter number of $W$ and $r$ is a positive integer  coprime with $h$ \cite{GG}. It appears that in general the corresponding parameter values have to be $c=(d-1)/h+m$, where $d$ is one of the degrees of $W$  and $m \in \Z_{\ge 0}$.  We explain how to construct all the singular  polynomials in the isotypic component of $V$ in terms of the Saito polynomials that are flat coordinates of the Saito metric. We use theory of Frobenius manifolds and especially Dubrovin's almost duality~\cite{Dubrovin_dual}. We show that singular polynomials under consideration correspond to the $W$-invariant polynomial twisted periods of the Frobenius manifold $V/W$, and we determine all such twisted periods.

 Firstly we prove that the first order derivatives of the Saito polynomials are singular polynomials at appropriate values of the parameter $c=(d-1)/h$ (Corollary~\ref{corSaitoPol}). Then we explain how to construct further singular polynomials with parameter $c$ shifted by an integer (Theorem~\ref{mainth}).
Then we show  in Corollary~\ref{all-sing}
that this construction provides all the singular polynomials in the isotypic component of the reflection representation. 

In the last section we present residue formulas for all the polynomial invariant twisted periods in the case of classical Coxeter groups $W$. Then we generalize them to get some singular polynomials for the complex reflection group  $W=S_n\ltimes(\mathbb Z/\ell\mathbb Z)^n$.

\section{Frobenius structures on the orbit spaces}

Let $V=\C^n$ with the standard constant metric $g$ given by $g(e_i,e_j)=(e_i,e_j)=\delta_{ij}$ where $e_i$, $i=1,\ldots,n$, is the standard basis in $V$. Let $(x^1,\ldots, x^n)$ be the corresponding orthogonal coordinates. Let $W$ be an irreducible finite Coxeter group of rank $n$ which acts in $V$ by orthogonal transformations such that $V$ is the complexified reflection representation of $W$. Let $\mathcal R \subset V$ be the Coxeter root system with the group $W$~\cite{Humphreys}. Let $y^1(x),\ldots, y^n(x)$ be a homogeneous  basis in the ring of invariant polynomials $S(V^*)^W=\C[x^1,\ldots,x^n]^W=\C[x]^W$. Let $d_\alpha$ be the corresponding degrees $d_\alpha = \deg y^\alpha(x)$, $\alpha=1,\ldots,n$. We assume the polynomials are ordered so that $d_1\ge\ldots\ge d_n$; $d_1=h$ is the Coxeter number of the group $W$. The polynomials $y^1,\ldots,y^n$ are coordinates on the orbit space ${\mathcal M}=V/W$. The Euclidean coordinates $x^1,\ldots, x^n$ can also be viewed as local coordinates on $\mathcal M\backslash\Sigma$, where $\Sigma$ is the discriminant set. Denote by $\Sympl=\{x\in V | (\g,x)=0 \mbox{ for some } \g \in {\mathcal R} \}$  the preimage of $\Sigma$ in the space $V$.

The metric $g$ is defined on ${\mathcal M}\setminus \Sigma$ due to its $W$-invariance. Let $g^{\alpha\beta}$ be the corresponding contravariant metric. Consider its Lie derivative
$\eta^{\alpha\beta}(y) = \partial_{y^1} g^{\alpha\beta}(y)$.
The metric $\eta$ is called the {\it Saito metric}. It is correctly defined (up to proportionality), and it is flat. There exist homogeneous coordinates $t^\a \in \C[x]^W$ ($1 \le \a \le n$), $\deg t^\a = d_\a$, such that $\eta$ is constant and anti-diagonal:
\begin{align*}
 &\eta^{\alpha\beta}=\partial_{t^1}g^{\alpha\beta}(t)=\delta^{\alpha+\beta}_{n+1}, &&1 \le \a,\b \le n,
\end{align*}
where $\delta^i_j=\delta_{ij}$ is the Kronecker symbol\footnote{We distinguish between the upper and lower indices as we will use the standard differential-geometrical convention that summation over the repeated upper and lower indices is assumed.}.
 Such coordinates are called {\it Saito polynomials}.

The pair of metrics $g, \eta$ forms a flat pencil which defines a Frobenius manifold~\cite{Dubrovin_lect}. We will be mainly concerned with the almost dual Frobenius structure~\cite{Dubrovin_dual}. It is defined by the prepotential
$$
F(x)= \frac{1}{2} \sum_{\g \in \mathcal R_+} (\g,x)^2 \log(\g,x),
$$
where summation is over the set of positive roots, and the roots are normalized so that $(\g,\g)=2$. The prepotential is quasi-homogeneous that is the Lie derivative
$$
{\mathcal L}_E F = \frac{2}{h} F + \mbox{quadratic terms in}  \,\,  x,
$$
where $E$ is the Euler vector field
$$
 E= \frac{1}{h} x^i\frac{\partial}{\partial x^i}= E^\alpha \frac{\partial}{\partial t^\alpha}
 $$
 with $E^\alpha = \frac{d_\alpha}{h}t^\alpha$.

Define the tensor
\begin{align}
 \Cstar_{i j k} = \frac{\partial^3 F}{\partial x^i \partial x^j\partial x^k}=\sum_{\g\in\mathcal R_+}\frac{\g_i\g_j\g_k}{(\g,x)}, \label{Cstar}
\end{align}
where $x \in V$ and $\gamma_i=(\gamma,e_i)$. Let $\Cstar^i_{jk}=g^{i l} \Cstar_{jkl}$. Then for any $x\in V\setminus \Sympl$ the tensor $\Cstar^i_{jk}=\Cstar^i_{jk}(x)$ forms the structure constants of an associative $n$-dimensional algebra~\cite{Dubrovin_dual}.

Define another tensor $C^{ij}_k= \Cstar^{ij}_k=g^{il} \Cstar^{j}_{kl}$ and consider the corresponding tensor with two low indices $C^\alpha_{\beta\lambda} = \eta_{\lambda\varepsilon} C^{\alpha\varepsilon}_\beta$ in the flat coordinates $t^\alpha$. The tensor $C^{\alpha}_{\beta\lambda}$ defines associative multiplication of tangent vectors at any point of the orbit space; the vector field $\p_{t^1}$ is the unity of this multiplication.
Let $C_\alpha$ be the corresponding $n\times n$ matrix with the entries $(C_\alpha)^\beta_\lambda = C^\beta_{\alpha\lambda}$. Let $U$ be the matrix  $U^\alpha_\beta=g^{\alpha\lambda}\eta_{\lambda\beta}$. The following result plays a key role.

\begin{theorem} \cite[Proposition 3.3]{Dubrovin_dual} \label{two-systems}
A function $p(t^1,\ldots,t^n)$ satisfies the system of equations
\begin{align}\label{second-derivatives-p}
\qquad\qquad \frac{\partial^2p}{\partial x^i \partial x^j} &= \nu\, \Cstar_{i j}^k \frac{\partial p}{\partial x^k}, &\qquad\qquad1\le i,j \le n 
\end{align}
if and only if the following equations hold
\begin{align}
 \xi_{\a}(t)&=\partial_{t^\alpha} p(t), &&1\le \a\le n,\label{xiDef}  \\
\partial_{t^\alpha}\xi(t)\, U
 &= \xi(t)\big(\nu +\Lambda\big) C_\a, &&1\le \a\le n, \label{ptt}
\end{align}
where $\xi(t)=\big(\xi_1(t),\ldots, \xi_n(t)\big)$ and $\Lambda$ is the diagonal matrix $$\Lambda=-\frac1{h}\Diag(d_1-1,\ldots,d_n-1).$$
\end{theorem}

Functions $p$ satisfying the system (\ref{second-derivatives-p}) are called {\it twisted periods} of the Frobenius manifold~\cite{Dubrovin_dual}.
We will be dealing with the system~\eqref{xiDef},~\eqref{ptt} so we note the following statement.

\begin{lemma} \label{lemFxiTp}
 For any polynomial solution $\xi(t)$ of the system~\eqref{ptt} there exists a polynomial $p(t)$ satisfying~\eqref{xiDef}.
\end{lemma}

\begin{proof}
 It is sufficient to check the compatibility $\partial_{t^\alpha}\xi_\beta=\partial_{t^\beta}\xi_\alpha$. Taking into account invertibility of the matrix $U$ on $\mathcal M\backslash\Sigma$ and using the equation~\eqref{ptt} we rearrange this equality as
\begin{align*}
\xi_\lambda\big(\nu +\Lambda\big)^\lambda_\varepsilon C^\varepsilon_{\beta\rho} (U^{-1})^\rho_\alpha=
\xi_\lambda\big(\nu +\Lambda\big)^\lambda_\varepsilon C^\varepsilon_{\alpha\rho} (U^{-1})^\rho_\beta.
\end{align*}
Note that it is necessary to check that
\begin{align} \label{CUinvCUinv}
 C^\varepsilon_{\beta\rho} (U^{-1})^\rho_\alpha=C^\varepsilon_{\alpha\rho} (U^{-1})^\rho_\beta.
\end{align}
The equality~\eqref{CUinvCUinv} is equivalent to
\begin{align} \label{CUCUeq}
 C^\varepsilon_{\lambda\rho} U^\lambda_\delta=C^\varepsilon_{\lambda\delta} U^\lambda_\rho.
\end{align}
Indeed, multiplying~\eqref{CUCUeq} by $(U^{-1})^\rho_\alpha (U^{-1})^\delta_\beta$ and summing by repeating indexes we obtain~\eqref{CUinvCUinv}. The relation~\eqref{CUCUeq}, in its turn, follows from the property $U^\a_\b = E^\varepsilon C^\a_{\varepsilon \beta}$ (see~\cite{Dubrovin_lect}) and the associativity conditions $C^\varepsilon_{\lambda\rho} C^\lambda_{\alpha\delta}=C^\varepsilon_{\lambda\delta} C^\lambda_{\alpha\rho}$.
\end{proof}

We observe that the Saito polynomials themselves satisfy the equations from Theorem~\ref{two-systems}. More exactly we have the following statement.

\begin{lemma} \label{lemptbeta}
For each $\beta=1,\ldots,n$ the function $p(t)=t^\beta$ is a solution to the system~\eqref{xiDef},~\eqref{ptt} with the parameter
\begin{align*}
 \nu=\frac{d_\beta-1}{h}.
\end{align*}
\end{lemma}

\begin{proof}
By substituting $\xi_\alpha=\delta^\beta_\alpha$ into~\eqref{ptt} we obtain the equations $0=e^\beta( \nu+\Lambda)C_\alpha$, where $e^\beta=(0,0,\ldots,1,\ldots,0)$ with $1$ in the site $\beta$. These equations are satisfied if $\nu-\frac{d_\beta-1}{h}=0$.
\end{proof}

Recall now the Dunkl operators associated with the Coxeter group $W$. We fix a parameter $c \in \C$. The Dunkl operator in the direction $e_i$ ($i=1,\ldots,n$) is given by
\begin{align}\label{dunkl-operator}
 \nabla_i=\partial_{x^i}-c\sum_{\gamma\in\mathcal R_+}\frac{\gamma_i}{(\gamma,x)}(1-s_\gamma),
\end{align}
where $s_\gamma$ denotes the orthogonal reflection with respect to the hyperplane $(\gamma,x)=0$. 
The key property of  Dunkl operators is their commutativity~\cite{Dunkl}:
\begin{align*}
 [\nabla_i,\nabla_j]&=0, &&1\le i,j\le n.
\end{align*}

\begin{proposition}\label{thm1}
 Suppose a $W$-invariant polynomial $p(x^1,\ldots,x^n)$ satisfies the system~\eqref{second-derivatives-p}. Then for any $j=1,\ldots,n$ the polynomial $v_j(x)=\partial_{x^j}p(x)$ satisfies the equations
\begin{align}
\nabla_i v_j&=0, &&i=1,\ldots,n, \label{Dbva}
\end{align}
where $\nabla_i$ is the Dunkl operator~\eqref{dunkl-operator} with the parameter $c=\nu$.
\end{proposition}

\begin{proof}
 By using the $W$-invariance of $p(x)$ we rearrange the left-hand side of the equation~\eqref{Dbva} as
\begin{equation*}
\nabla_i v_j =\partial_{x^i}\partial_{x^j} p(x)-\nu\sum_{\gamma\in\mathcal R_+}\frac{\gamma_i}{(\gamma,x)}\big(\partial_{x^j} p(x)-\frac{\partial}{\partial(s_\gamma e_j)}p(x)\big)=
\end{equation*}
\begin{equation*}=\partial_{x^i}\partial_{x^j} p(x)-\nu\sum_{\gamma\in\mathcal R_+}\frac{2\gamma_j\gamma_i}{(\gamma,\gamma)(\gamma,x)}\Big(\gamma,\frac{\partial}{\partial x}\Big)p(x).
\end{equation*}
By using $\big(\gamma,\frac{\partial}{\partial x}\big)=\sum_i\gamma_i\partial_{x^i}$ and the formula~\eqref{Cstar} we obtain
\begin{align}
\nabla_i v_j&=\partial_{x^i}\partial_{x^j} p(x)-\nu\Cstar^k_{ij}\partial_{x^k} p(x).
\end{align}
Thus the property~\eqref{Dbva} follows from~\eqref{second-derivatives-p}.
\end{proof}

\begin{corollary} \label{corSaitoPol}
 Consider the Saito polynomial $t^\beta=t^\beta(x)$ for some $\beta=1,\ldots,n$.  Then the derivatives $v_j=\partial_{x^j} t^\beta(x)$ satisfy the relations~\eqref{Dbva}, that is
\begin{align}
 \nabla_i\nabla_j t^\beta(x)=\nabla_i\partial_{x^j} t^\beta(x)&=0, &&i,j=1,\ldots,n,\label{Dbatbeta}
\end{align}
if the Dunkl operators have parameter $c=\frac{d_\beta-1}{h}$.
\end{corollary}

\begin{definition} \cite{DJO}
A polynomial $q(x)$ is called {\it singular} if $\nabla_i q(x)=0$ for all $i=1,\ldots,n$.
\end{definition}
Thus Corollary~\ref{corSaitoPol} deals with the singular polynomials $v_j$. The $W$-module  $\langle v_1,\ldots, v_n\rangle$ is isomorphic to the reflection representation of the Coxeter group $W$.

\section{Shifting}

In the previous Section we established that derivatives of the Saito polynomials $t^\beta$ are singular polynomials for the appropriate values of the parameter $c=c_\beta=(d_\beta-1)/h$. In this Section we explain how to generate further singular polynomials starting with Saito polynomials. The corresponding parameters $c$ differ from $c_\beta$ by integers.

We start with a known property of the solutions of the system~\eqref{ptt}.

\begin{lemma} \cite[Lemma 3.6]{Dubrovin_dual} \label{lempt1}
 If $\xi(t)$ is a solution of the system~\eqref{ptt} then $\wt\xi(t)=\partial_{t^1}\xi(t)$ is a solution of the same system with $\nu$ replaced by $\nu-1$:
\begin{align*}
\partial_{t^\alpha}\wt\xi(t)\, U
 &=\wt\xi(t)\, \big(\nu-1+ \Lambda\big) C_\alpha, &&1\le \alpha\le n. 
\end{align*}
\end{lemma}

Note also that  if a function $p(t)$ is a solution of the system~\eqref{xiDef},~\eqref{ptt} then $\partial_{t^1} p(t)$ is a solution of the same system with $\nu$ replaced by $\nu-1$. Indeed, the partial derivatives of the function $\partial_{t^1} p(t)$ are $\partial_{t^\alpha}\partial_{t^1} p(t)=\partial_{t^1} \xi_\alpha(t)$, so they satisfy  the system~\eqref{ptt} with $\nu-1$.

The new solution $\wt\xi(t)$ given in Lemma~\ref{lempt1} can be represented as $\wt\xi=\xi \big(\nu+\Lambda\big)U^{-1}$. If $\nu\ne\dfrac{d_\alpha-1}h$ for all $\alpha=1,\ldots,n$ then the matrix $\nu+ \Lambda$ is invertible and we can rewrite this relation as $\xi=\wt\xi\, U\big(\nu+\Lambda\big)^{-1}$. This suggests how to invert the Lemma~\ref{lempt1}  in order to generate solutions with increased value of $\nu$.

\begin{lemma} \label{lempt1conv}
Let $\xi(t)$ be a solution of the system~\eqref{ptt}. Assume that $\nu\ne\dfrac{d_\alpha-1}h-1$ for all $\alpha=1,\ldots,n$. Then
\begin{align} \label{whxixi}
 \wh\xi(t)=\xi(t) U\big(\nu+1+ \Lambda\big)^{-1}
\end{align}
is a solution of the system~\eqref{ptt} with $\nu$ replaced by $\nu+1$:
\begin{align}
\partial_{t^\alpha}\wh\xi(t)\, U
 &=\wh\xi(t)\, \big(\nu+1+ \Lambda\big) C_\alpha, &&1\le \alpha\le n. \label{pttHat}
\end{align}
\end{lemma}

\begin{proof}
 Let $t_0$ be a generic point in $\mathcal M$ and let  $\xi(t_0)$ be the value of $\xi(t)$ at this point. Then the value of $\wh\xi(t)$ at this point is $\wh\xi(t_0)=\xi(t_0) U(t_0)\big(\nu+1+\Lambda \big)^{-1}$. There exists a solution $\wh\zeta(t)$ of the system~\eqref{pttHat} such that $\wh\zeta(t_0)=\wh\xi(t_0)$. By Lemma~\ref{lempt1} the covector $\zeta(t)=\wh\zeta(t)\big(\nu+1+\Lambda\big)U^{-1}$ is a solution of~\eqref{ptt}. Note that there exists a unique solution of the system~\eqref{ptt} with a given value in the point $t_0$. So taking into account $\zeta(t_0)=\wh\xi(t_0)\big(\nu+1+\Lambda\big)U^{-1}(t_0)=\xi(t_0)$ one yields $\zeta(t)=\xi(t)$. Therefore $\widehat\zeta(t)=\zeta(t)U\big(\nu+1+\Lambda\big)^{-1}=\xi(t)U\big(\nu+1+\Lambda\big)^{-1}=\wh\xi(t)$ and $\wh\xi(t)$ satisfies the equation~\eqref{pttHat}.
\end{proof}

\begin{remark}\label{rem-increase-p}
Suppose $\xi(t)$ in Lemma~\ref{lempt1conv} is polynomial. Then by Lemma~\ref{lemFxiTp}
there exists a polynomial $\wh p(t)$ such that $\wh\xi_\alpha(t)=\partial_{t^\alpha}\wh p(t)$. Thus the covector $\wh\xi(t)$ satisfies the whole system~\eqref{xiDef},~\eqref{ptt} (with $\nu$ replaced by $\nu+1$)  with some polynomial $\wh p(t)$.
\end{remark}

By applying Lemma~\ref{lempt1conv} to the first order derivatives of the Saito polynomials we get the following result.

\begin{proposition} \label{propwhxibeta}
 The covector $\wh\xi$ with the components
\begin{align} \label{whxibeta}
 \wh\xi_\alpha=\frac{U^\beta_{\alpha}}{d_\beta-d_\alpha+h},   \,\,\,\, \a=1,\ldots,n,
\end{align}
satisfies the equations
\begin{equation}\label{system-hat}
\partial_{t^\alpha}\widehat\xi(t)\, U
 = \widehat\xi(t)\big(\widehat\nu +\Lambda\big) C_\alpha, \,\,\,\,\, 1\le \a\le n,
\end{equation}
with $\widehat\nu=\dfrac{d_\beta-1}h+1$.
\end{proposition}

\begin{proof}
 Let $\nu=\dfrac{d_\beta-1}h$ and consider the solution $p(t)=t^\beta$ of the system~\eqref{xiDef},~\eqref{ptt} and the corresponding covector $\xi$ with the components $\xi_\alpha=\delta^\beta_\alpha$ (see Lemma~\ref{lemptbeta}).
 By Lemma~\ref{lempt1conv} the  covector $\wh\xi$ given by the formula~\eqref{whxixi} is a solution of~\eqref{system-hat} for $\wh\nu=\nu+1$. By substituting $\xi_\alpha=\delta^\beta_\alpha$ into~\eqref{whxixi} we obtain the components~\eqref{whxibeta}.
\end{proof}

This leads us to the following result.
\begin{theorem}
For any $\zeta \in V$, $\beta=1,\ldots,n$
the polynomial
\begin{align} \label{qtheor}
q(x)=\sum_{i,\a=1}^n  \frac{1}{d_\beta - d_\a+h} \frac{\partial t^\beta}{\partial x^i} \frac{\partial t^{n+1-\a}}{\partial x^i}{\partial_\zeta t^\a}
\end{align}
is a singular polynomial for the Dunkl operators with the parameter $c=\dfrac{d_\beta-1}h+1$.
\end{theorem}

\begin{proof}
First we rearrange
\begin{align} \label{U-expand}
U^\beta_\alpha = g^{\beta \lambda} \eta_{\lambda \alpha}=g^{\beta, n+1-\alpha}=
\sum_{a=1}^n\frac{\partial t^\beta}{\partial x^a}\frac{\partial t^{n+1-\alpha}}{\partial x^a}.
\end{align}
It follows from Proposition~\ref{propwhxibeta} and Lemma~\ref{lemFxiTp} that there exists a $W$-invariant polynomial $p(x)$ such that
$$
\partial_{t^\alpha} p =  \frac{1}{d_\beta-d_\alpha+h}
\sum_{a=1}^n \frac{\partial t^\beta}{\partial x^a}\frac{\partial t^{n+1-\alpha}}{\partial x^a}.
$$
By Proposition~\ref{thm1} the derivative $q(x)=\partial_\zeta p(x)$ is a singular polynomial. It has the required form as $\partial_\zeta p = \partial_{t^\a}p  \partial_{\zeta} t^\a$.
\end{proof}

As an example consider the case $\beta=n$. The corresponding Saito polynomial $t^n$ is proportional to $(x^{1})^2+(x^{2})^2+\ldots+(x^{n})^2$. The right-hand side of the equality~\eqref{qtheor} is then proportional to
$$
\sum_{\a=1}^n t^{n+1-\a} \partial_\zeta t^\a
$$
as the polynomial $t^{n+1-\a}$ is homogeneous of degree $d_{n+1-\a}$ and $d_\a+d_{n+1-\a}=h+2$.
We arrive at the following corollary.
\begin{proposition}\label{prophhp1}
For any $\zeta \in V$
the polynomial
\begin{align} \label{qcor}
q(x)=\partial_\zeta \sum_{\alpha=1}^n t^\alpha t^{n+1-\alpha}
\end{align}
is a singular polynomial for the Dunkl operators with the parameter $c=(h+1)/h$.
\end{proposition}

\begin{remark}
For $c=1/h+m, m \in \Z_{\ge 0}$ the existence of the singular polynomials in the isotypic component of the  reflection representation is known from~\cite{BEG} (see also~\cite{Go}). Further, in the case of Weyl groups $W$ and $c=r/h$ where $r$ is a positive integer coprime with $h$,  the existence of singular polynomials in the isotypic component of the reflection representation is known from \cite{GG}.
\end{remark}

\begin{example}
Let ${\mathcal R}={\mathcal A}_{n} \subset V \subset \C^{n+1}$ be given in its standard embedding so that $V\subset \C^{n+1}$ is defined by $\sum_{i=1}^{n+1} z_i=0$ where $z_i$ are standard coordinates in $\C^{n+1}$. Define $s^\alpha = \mathop{\mathrm{Res}}\nolimits_{z=\infty} \prod_{i=1}^{n+1} (z-z_i)^\frac{n+1-\alpha}{n+1} d z$ for $\alpha=1,\ldots,n$. Then Saito coordinate $t^\alpha = \frac{1}{n-\a+1}s^\alpha|_V$ (see~\cite{Dubrovin_lect},~\cite{SYS}). The polynomials $s^\alpha$ satisfy  $\sum_{i=1}^{n+1} \frac{\partial s^\alpha}{\partial z_i} =0$ so for any $\zeta \in \C^{n+1}$, $i=1,\ldots,n+1$ Corollary~\ref{corSaitoPol} gives
$\nabla^{\frac{n+1-\alpha}{n+1}}_i \partial_\zeta s^\alpha=0$,
 where
$$
\nabla^c_i =  \frac{\partial}{\partial z_i} - c\sum_{\nad{j=1}{j \ne i}}^{n+1} \frac{1-s_{ij}}{z_i-z_j}
$$
with $s_{ij}$ exchanging $z_i$ and $z_j$ (see also~\cite[Proposition 11.14]{Et},~\cite{Du} where this fact was established). Further, Proposition~\ref{prophhp1} gives that
the polynomial $q(z_1,\ldots,z_{n+1})=\partial_\zeta \sum_{\alpha=1}^n t^\alpha t^{n+1-\alpha}$ satisfies $\nabla_i ^{\frac{h+1}{h}}q =0$.
\end{example}

By iterating the previous arguments we get the following statement.

\begin{theorem}\label{mainth}
Let $m\in \Z_{\ge0}$, fix $\beta$,  $1\le \beta \le n$. Define the covector $\xi^{(m)}=(\xi^{(m)}_1,\ldots,\xi^{(m)}_n)$ by
\begin{align} \label{xim}
 \xi^{(m)}= \xi^{(0)} \mathop{\overrightarrow\prod}_{1\le j\le m} U\, \big( \Lambda + \frac{d_\beta-1}{h}+j \big)^{-1},
\end{align}
where $\xi^{(0)}$ has components $\xi^{(0)}_\alpha = \delta^{\beta}_\alpha$, $\alpha=1,\ldots,n$ and the factors are ordered as $\mathop{\overrightarrow\prod}\limits_{1\le j\le m}A_j=A_1A_2\cdots A_m$.
Then for any $i=1,\ldots,n$ the polynomials
\begin{align} \label{qi}
q_i(x) =q_{\beta,i}(x) =
\sum_{\a=1}^n \xi^{(m)}_\a \frac{\partial t^\a}{\partial x^i}
\end{align}
are singular polynomials for the Dunkl operators with the parameter $c=\frac{d_\beta-1}{h}+m$. These polynomials are homogeneous of degree $d_\b-1+h m$.
\end{theorem}

\section{Singular polynomials in the reflection representation}

We are going to show that polynomials (\ref{qi}) generate all singular polynomials in the isotypic component of the reflection representation of $W$. Firstly we note that each copy of the reflection representation spanned by the singular polynomials is governed by a single $W$-invariant polynomial.

\begin{proposition}\label{Qexists}
Let a subspace $M_0\subset\mathbb C[x]$ be spanned by the singular polynomials and suppose that $M_0\cong V$ as $W$-modules. Choose a basis $\{P_1,\ldots, P_n\}$ in $M_0$ such that each polynomial $P_i$ is mapped to the basis vector $e_i \in V$ under the isomorphism. Then there exists $Q \in \C[x^1,\ldots,x^n]^W$ such that $\frac{\partial Q}{\partial x^i}= P_i$ for all $i=1,\ldots, n$.
\end{proposition}
\begin{proof}
We have $\nabla_i P_j = \nabla_j P_i =0$ for all $i,j=1,\ldots,n$. Hence
\begin{align} \label{crosspartial}
\partial_{x^i} P_j - \partial_{x^j} P_i - c \sum_{\gamma\in {\mathcal R}_+} \frac{(\gamma,e_i)}{(\gamma,x)}(1-s_\gamma)P_j +  c \sum_{\gamma\in {\mathcal R}_+} \frac{(\gamma,e_j)}{(\gamma,x)}(1-s_\gamma)P_i=0.
\end{align}
Since $(1-s_\gamma)e_i = \frac{2(\gamma,e_i)}{(\gamma,\gamma)}\gamma$ we get
$$
(\gamma,e_i)(1-s_\gamma)P_j = (\gamma,e_j)(1-s_\gamma)P_i
$$
for any $\gamma \in {\mathcal R}_+$.
Thus it follows from the relation~\eqref{crosspartial} that  $\partial_{x^i} P_j = \partial_{x^j} P_i$ so there exists $Q\in \C[x^1,\ldots,x^n]$ such that $\frac{\partial Q}{\partial x^i}=P_i$ for all $i=1,\ldots,n$. Let us also check that $Q$ is $W$-invariant. Fix $\gamma\in {\mathcal R}_+$. Then for any $i=1,\ldots,n$ we have $s_\gamma P_i = \partial_{s_\gamma e_i} (s_\gamma Q)$, and on the other hand $s_\gamma P_i = \partial_{s_\gamma e_i} Q$. Thus $\partial_{s_\gamma e_i}(Q-s_\gamma Q)=0$, so $Q=s_\gamma Q$ as required.
\end{proof}

\begin{corollary}\label{corSaitoPolgeneral}
The singular polynomials~\eqref{qi} can be represented as
\begin{align} \label{qiqi}
q_i = \frac{\partial Q}{\partial x^i},
\end{align}
where
\begin{align} \label{QQ}
Q=Q_\beta = \frac{1}{d_\b+h m}
\sum_{\a=1}^n d_\a \xi^{(m)}_\a  t^\a,
\end{align}
$\,i=1,\ldots,n$, and we keep notations of Theorem \ref{mainth}.
\end{corollary}
\begin{proof}
By Proposition~\ref{Qexists} we have the relation~\eqref{qiqi} for some invariant polynomial $Q$ of degree $d_\b +h m$. Hence $(d_\b+h m) Q=\sum_{i=1}^n x^i q_i$
\end{proof}

\begin{remark}
It is recently explained in~\cite{KL} how ${\mathcal N}=4$ multi-particle mechanical system with $D(2,1;\a)$ superconformal symmetry can be constructed based on a solution of the WDVV equations and a particular twisted period. It follows from Theorem~\ref{mainth} and Corollary~\ref{corSaitoPolgeneral} that the polynomials $q_i, Q$  given by~\eqref{qi},~\eqref{QQ} define a superconformal mechanical system with the bosonic potential proportional to
$
Q^{-2} \sum_{i=1}^n q_i^2
$
at the parameter value $\a=-\frac12 (d_\b + h m)$.
\end{remark}

\begin{remark}
Let $g(x^1,\ldots,x^n)$ be a homogeneous  $W$-invariant polynomial of positive degree. Let $L_g$ be the differential operator which acts on the $W$-invariant functions by $g(\nabla_1,\ldots,\nabla_n)$. The operators $L_g$ commute with each other and include the corresponding {\it Calogero-Moser operator}~\cite{Et}.  It follows from Corollary~\ref{corSaitoPolgeneral} that if $c=\frac{d_\beta-1}{h}+m$ then  $L_g Q_\beta=0$ so the polynomial $Q_\beta$ belongs to the joint kernel of the Calogero-Moser operator and its quantum integrals. In particular, the Saito polynomial $t^\beta$ satisfies $L_g t^\beta=0$ if $c=\frac{d_\beta-1}{h}$.
\end{remark}

Now we move to the main statement of this Section on possible polynomial twisted periods of the Frobenius manifold $\mathcal M$.

\begin{theorem}\label{mthsect4uniqueness}
Let $L$ be the linear space of solutions $p(x)$ to the system~\eqref{second-derivatives-p} such that $p \in \C[x]^W$. Then $\dim L =1$ unless
$\nu=\dfrac{d_\beta-1}{h}+m$ for some degree $d_\beta$, and $m\in\mathbb Z_{\ge0}$. In this case $\dim L=2$ unless $W=D_n$ where $n$ is even and $d_\b=n$ when $\dim L=3$.
\end{theorem}

\begin{proof}
Suppose $p\in \C[x]^W$ is a homogeneous solution of the system~\eqref{second-derivatives-p} such that $D=\deg p >0$. By Proposition~\ref{thm1} the polynomials $v_i(x)= \frac{\partial p(x)}{\partial x^i}$ are singular at the parameter $c=\nu$. It follows from the relation
\begin{align*}
 \sum_{i=1}^n x^i\nabla_i=\sum_{i=1}^n x^i\partial_{x^i}-\nu\sum_{\gamma\in \mathcal R_+}(1-s_\gamma)
\end{align*}
that $\nu= \frac{D-1}{h}\ge0$
(c.f.~\cite{DJO}).
Consider firstly the case $0\le\nu\le1$.
Equation~\eqref{ptt} at $\a=1$ takes the form $0=\xi(t)(\nu+\Lambda)$ since $C_1=Id$ and $\deg\xi_\lambda=D-d_\lambda<h$ for any $1 \le \lambda \le n$.
Hence the matrix $\nu+\Lambda$ is degenerate
so $\nu= \frac{d_\beta-1}{h}$ for some $\beta$, and $D=d_\b$. Moreover $p(x)$ as a polynomial of the Saito coordinates has to be a linear combination of $t^\lambda$ with $d_\lambda=d_\b$.

Now let $\nu>1$. The polynomial $\p_{t^1}p(t)$ is non-constant as the matrix $\nu+\Lambda$ is non-degenerate. By Lemma~\ref{lempt1} the polynomial $\p_{t^1}p(t)$ is a solution of the system~\eqref{second-derivatives-p} with $\nu$ replaced by $\nu-1$. It follows by Lemma~\ref{lempt1conv} and Remark~\ref{rem-increase-p} that the spaces of positive degree homogeneous invariant solutions of~\eqref{second-derivatives-p} for the parameter $\nu$ and $\nu$ replaced with $\nu-1$ are isomorphic.
\end{proof}

The arguments from the proof of Proposition~\ref{thm1} can be reversed so it follows that the polynomial~\eqref{QQ} is a twisted period with $\nu=(d_\b-1+h m)/h$ of the Frobenius manifold $\mathcal M$. Thus Theorem~\ref{mthsect4uniqueness} implies the following statement.

\begin{corollary} \label{all-twist}
Let $p \in \C[x]^W$ be a non-constant  twisted period of the Frobenius manifold $\mathcal M$. Then $\nu =(d_\b-1+h m)/h$ for some degree $d_\b$ and $m \in \Z_{\ge 0}$. Further, $p= \lambda Q_{\beta}$ unless $W=D_{2r}$ and $d_\b=2r$ for some $r \in \N$ in which case $p=\lambda Q_{r}+\mu Q_{r+1}$, where $Q_{\beta}$ are given by~\eqref{QQ} and $\lambda,\mu\in \C$.
\end{corollary}

Consider now a linear space $M_0$ of singular polynomials such that $M_0 \cong V$ as $W$-modules. Then $Q\in\mathbb C[x]^W$ defined by Proposition~\ref{Qexists} is a twisted period and hence Corollary~\ref{all-twist} allows to describe all such $W$-modules $M_0$.

\begin{corollary} \label{all-sing}
Let $q$ be a homogeneous singular polynomial. Suppose the linear space spanned by the polynomials $wq$, $w\in W$, is isomorphic to $V$ as $W$-module. Then $\deg q = d_\b -1 + h m$ for some degree $d_\b$ of $W$,  $m \in \Z_{\ge 0}$, and $c=(d_\b-1+h m)/h$.

Further, $q=\sum_{i=1}^n \eta_i q_{\beta,i}$ unless $W=D_{2r}$ and $d_\b=2r$ for some $r \in \N$ in which case $q=\sum_{i=1}^n  \eta_i (\lambda q_{r,i}+ \mu  q_{r+1,i})$, where $q_{\beta,i}$ are given by~\eqref{qi},~\eqref{xim} and $\lambda,\mu, \eta_i \in \C$.

In the former cases all the homogeneous singular polynomials for $c=(d_\b-1+h m)/h$ in the isotypic component of the reflection representation are described as linear combinations  $\sum_{i=1}^n \eta_i q_{\beta,i}$ while in the latter case the homogeneous singular polynomials in the isotypic component of the reflection representation form the $2n$-dimensional subspace of polynomials $\sum_{i=1}^n \lambda_i q_{r,i}+\mu_i q_{r+1,i}$, where $\lambda_i, \mu_i \in \C$.
\end{corollary}

\section{Further examples for classical series}

While we express singular polynomials in the isotypic component of the reflection representation of $W$ through the Saito polynomials, in certain cases direct formulas exist. We refer to~\cite{DJO},~\cite{Chm} for the case of dihedral groups, and to~\cite{Du},~\cite{EtChm} for the case when $W$ is of classical type.
For instance, it follows from~\cite{Du} and Corollary \ref{all-twist} that {\it all} the polynomial invariant twisted periods for $W=A_{n}$ are given by
$$
Q=\mathop{\mathrm{Res}}\nolimits_{z=\infty}  \prod_{j=1}^{n+1} (z-z_j)^\nu dz  \Big|_{\sum z_j=0},
$$
where $\nu=\frac{s}{n+1}+m$ with $s=1,\ldots,n$ and $m \in \Z_{\ge 0}$ (c.f.~\cite{Dubrovin_dual}).
Further, let
\begin{equation} \label{twistedperiods}
Q=\mathop{\mathrm{Res}}\nolimits_{z=\infty} z^{a}\prod_{j=1}^n\big(z^2-x_j^2\big)^\nu dz,
\end{equation}
where $x_1,\ldots,x_n$ are standard coordinates in $\C^n$.
Then {\it all} the polynomial invariant twisted periods for $W=B_n$ have the form~\eqref{twistedperiods} where $a=0$, $\nu=\frac{2 s -1}{2n}+m$ with $s=1,\ldots,n$ and $m \in \Z_{\ge 0}$.
Similar, \eqref{twistedperiods}~is the twisted period for $W=D_n$ if $a=-2\nu$, $\nu=\frac{2 s -1}{2(n-1)}+m$ with $s=1,\ldots,n-1$ and $m \in \Z_{\ge 0}$. {\it All} the remaining polynomial invariant twisted periods for $W=D_n$ have the form
\begin{equation} \label{twistedperiodsDrest}
Q=\mathop{\mathrm{Res}}\nolimits_{z=0} z^{-2m-1}\prod_{j=1}^n\big(z^2-x_j^2\big)^{m+\frac12} dz,
\end{equation}
where $m \in \Z_{\ge 0}$ (c.f~\cite{Egouchi}).
We note that for even $n$ the polynomial $Q=Q_{m,\infty}$ given by~\eqref{twistedperiods} with $a=-2m-1$, $\nu = m+\frac12$, $m\in \Z_{\ge 0}$ has same degree as the polynomial
 $Q=Q_{m,0}$ given by~\eqref{twistedperiodsDrest} with the same $m$. These polynomials are not proportional as for $m>0$ the polynomial $\p_{t^1} Q_{m,\infty}$ is a nonzero multiple of the polynomial $Q_{m-1,\infty}$, and  $\p_{t^1} Q_{m,0}$ is a nonzero multiple of $Q_{m-1,0}$.

This leads to the following proposition which can also be checked directly.

\begin{proposition} \label{propnum12}
Let polynomial $Q(x)$ be given by~\eqref{twistedperiodsDrest}
with  $m\in \Z_{\ge 0}$. For any $\zeta \in V$ the polynomial $\p_\zeta Q$ is singular for $W=D_n$ with the parameter $c=m+\frac12$.
\end{proposition}

Proposition~\ref{propnum12} can be generalized to the case of complex reflection group $W=S_n\ltimes(\mathbb Z/\ell\mathbb Z)^n$, where $\ell \in \Z$ satisfies $\ell \ge 2$. This group is generated by the reflections $\sigma_{ij}^{(a)}$ and $s_i$  acting on the standard basis as
\begin{align}
 \sigma_{ij}^{(a)}e_i&=\omega^{-a}e_j, &\sigma_{ij}^{(a)}e_j&=\omega^{a}e_i,
 &\sigma_{ij}^{(a)}e_k&=e_k, \\
 s_ie_i&=\omega e_i, &s_i e_k&=e_k,
\end{align}
where $\omega=e^{2\pi i/\ell}$ is the primitive root of unity, $i,j,k=1,\ldots,n$, $i\ne j\ne k\ne i$, $a=0,\ldots,\ell-1$.
The Dunkl operators in this case have the form~\cite{DO}
\begin{align}\label{dunklgmp}
 \nabla_i=\partial_{x_i}-\nu\sum_{a=0}^{\ell-1}\sum_{\nad{j=1}{j\ne i}}^n\frac1{x_i-\omega^ax_j}\big(1-\sigma_{ij}^{(a)}\big)-\sum_{b=1}^{\ell-1}c_b\sum_{a=0}^{\ell-1}\frac{\omega^{-ab}}{x_i}s_i^a,
\end{align}
where  $\nu,c_1,\ldots,c_{\ell-1}\in\mathbb C$ (we do not suppose any more that parameters of the Dunkl operators are equal). Denote $c_0=0$ and $c_{a\ell+b}=c_b$ for $a\in\mathbb Z$, $b=0,1,\ldots,\ell-1$.

The next statement generalizes Proposition~\ref{propnum12}. The form of the singular polynomials is suggested by~\cite{EtChm} where some singular polynomials for the group $W$ were found using the residues at infinity.

\begin{proposition}\label{complex-gropus}
 Let $q\in \Z$ satisfy $1\le q\le\ell-1$. Suppose that $\nu=m+\frac{\ell-q+s}\ell$, $c_{q-s}=0$ and $c_{-s}=\frac{s}\ell$ for some $m,s\in\mathbb Z_{\ge0}$. Then the formulae
\begin{align}
 f_j=\prod_{i=1}^nx_i^{\ell \nu}\mathop{\mathrm{Res}}\nolimits_{z=0} z^{-\ell m-1}\prod_{i=1}^n\bigg(1-\frac{z^\ell}{x_i^\ell}\bigg)^{\nu}\,\frac{x_j^{q}dz}{x_j^\ell-z^\ell}, \label{fi}
\end{align}
define singular polynomials that is $\nabla_i f_j=0$ for any $i, j=1,\ldots,n$, where the operator $\nabla_i$ is given by~\eqref{dunklgmp} (and it is assumed in~\eqref{fi} that $x_i\ne 0$ for all $i$). These polynomials are homogeneous of degree $(n-1)(m\ell+\ell-q)+ns$ and they span an irreducible $n$-dimensional representation of $W$.
\end{proposition}


\begin{proof}
By calculating the residue in~\eqref{fi} explicitly one yields
\begin{align*}
 f_j= (-1)^m\sum_{k_1,\ldots,k_n=0\atop k_1+\ldots+k_n=m}^\infty \binom{\nu-1}{k_j} x_j^{\ell(m-k_j)+s}\prod_{i=1\atop i\ne j}^n\binom{\nu}{k_i} x_i^{\ell(\nu-k_i)},
\end{align*}
where $\binom{\alpha}k=\alpha(\alpha-1)\cdots(\alpha-k+1)/k!$ and $\binom{\alpha}0=1$. Since $\nu>m$ the function $f_j$ is a homogeneous polynomial of $x_1,\ldots,x_n$, it has degree $(n-1)\ell \nu+s=(n-1)(m\ell+\ell-q)+ns$. Note that all the coefficients in this expression do not vanish and hence $f_j\ne0$. The generators of $W$ act on the polynomials~\eqref{fi} by the formulae
\begin{align}
 \sigma_{ij}^{(a)}f_i&=\omega^{aq}f_j, &\sigma_{ij}^{(a)}f_j&=\omega^{-aq}f_i,
 &\sigma_{ij}^{(a)}f_k&=f_k, \label{sigmaf} \\
 s_if_i&=\omega^{-s} f_i, &s_i f_k&=\omega^{q-s} f_k, \label{sif}
\end{align}
where $i,j,k=1,\ldots,n$, $i\ne j\ne k\ne i$, $a=0,\ldots,\ell-1$. It follows that the space spanned by $f_j$, $j=1,\ldots,n$, is $n$-dimensional irreducible representation of $W$.

Now let us show that the polynomial $f_j$ is singular. Let $i\ne j$, then using the formulae~\eqref{sigmaf},~\eqref{sif} we derive
\begin{align}
 \nabla_i f_j=\partial_{x_i}f_j-\nu\sum_{a=0}^{\ell-1}\frac1{x_i-\omega^ax_j}\big(f_j-\omega^{-aq}f_i\big)-\sum_{b=1}^{\ell-1}c_b\sum_{a=0}^{\ell-1}\frac{\omega^{-ab}}{x_i}\omega^{a(q-s)}f_j. \label{nablafj}
\end{align}
The last term in~\eqref{nablafj} equals $-\ell c_{q-s}x_i^{-1}f_j=0$ as $\sum_{a=0}^{\ell-1}\omega^{ab}=0$ for $b\notin\mathbb\ell\mathbb Z$. Thus $\nabla_i f_j=\prod_{i=1}^nx_i^{\ell \nu}\mathop{\mathrm{Res}}\nolimits_{z=0} z^{-\ell m-1}\prod_{i=1}^n\big(1-z^\ell/x_i^\ell\big)^{\nu}F_{ij}dz$
where
\begin{align}
 F_{ij}&=\nu\frac{\ell x_i^{\ell-1}x_j^q}{(x_i^\ell-z^\ell)(x_j^\ell-z^\ell)}-
\nu\sum_{a=0}^{\ell-1}\frac1{x_i-\omega^ax_j}\Big(\frac{x_j^q}{x_j^\ell-z^\ell}-\frac{\omega^{-aq}x_i^q}{x_i^\ell-z^\ell}\Big)= \notag \\
&=-\frac{\nu}{(x_i^\ell-z^\ell)(x_j^\ell-z^\ell)}\Big(-\ell x_i^{\ell-1}x_j^q+ \notag \\
&+\sum_{a=0}^{\ell-1}\frac{z^\ell\omega^{-aq}(x_i^q-\omega^{aq}x_j^q)+x_i^qx_j^q(x_i^{\ell-q}-\omega^{-aq}x_j^{\ell-q})}{x_i-\omega^ax_j}\Big)= \notag \\
&=-\frac{\nu}{(x_i^\ell-z^\ell)(x_j^\ell-z^\ell)}\Big(-\ell x_i^{\ell-1}x_j^q+ \notag \\
&+\sum_{a=0}^{\ell-1}z^\ell\omega^{-aq}\sum_{b=0}^{q-1}x_i^b\omega^{a(q-1-b)}x_j^{q-1-b}+
\sum_{a=0}^{\ell-1}x_i^qx_j^q\sum_{b=0}^{\ell-q-1}x_i^{\ell-q-1-b}\omega^{ab}x_j^b\Big). \notag 
\end{align}
Thus the first double sum in $F_{ij}$ vanishes and the last double sum equals $\ell x_i^{\ell-1}x_j^q$. Thus $F_{ij}=0$ and $\nabla_i f_j=0$ for $i\ne j$.

Now to prove that $\nabla_j f_j=0$ it is sufficient to check $\sum_{i=1}^nx_i\nabla_i f_j=0$. Since $f_j$ is homogeneous of order $(n-1)\ell \nu+s$ we have $\sum_{i=1}^nx_i\partial_{x_i} f_j=\big((n-1)\ell \nu+s\big)f_j$. Hence using $c_{q-s}=0$ as previously one yields
\begin{align}
 \sum_{i=1}^nx_i\nabla_i f_j=\big((n-1)\ell \nu+s\big)f_j+\prod_{k=1}^n x_k^{\ell \nu} \mathop{\mathrm{Res}}\nolimits_{z=0} z^{-\ell m-1}\prod_{k=1}^n(1-z^\ell/x_k^\ell)^{\nu}F_jdz, \notag
\end{align}
where
\begin{multline*}
 F_j=-\nu\sum_{a=0}^{\ell-1}\sum_{i=1\atop i\ne j}^n\frac{x_i}{x_i-\omega^ax_j}\Big(\frac{x_j^q}{x_j^\ell-z^\ell}-\frac{\omega^{-aq}x_i^q}{x_i^\ell-z^\ell}\Big)-\\
-\nu\sum_{a=0}^{\ell-1}\sum_{k=1\atop k\ne j}^n\frac{x_j}{x_j-\omega^ax_k}\Big(\frac{x_j^q}{x_j^\ell-z^\ell}-\frac{\omega^{aq}x_k^q}{x_k^\ell-z^\ell}\Big)-\sum_{b=1}^{\ell-1}c_b\sum_{a=0}^{\ell-1}\frac{\omega^{-ab-as}x_j^q}{x_j^\ell-z^\ell}.
\end{multline*}
By taking into account $\frac{x_i}{x_i-\omega^a x_j}+\frac{x_j}{x_j-\omega^{-a} x_i}=1$ we obtain
\begin{multline*}
 F_j=-\nu\sum_{a=0}^{\ell-1}\sum_{i=1\atop i\ne j}^n\Big(\frac{x_j^q}{x_j^\ell-z^\ell}-\frac{\omega^{-aq}x_i^q}{x_i^\ell-z^\ell}\Big)-\ell c_{-s}\frac{x_j^q}{x_j^\ell-z^\ell}= \\
=-\frac{x_j^q}{x_j^\ell-z^\ell}\big((n-1)\ell \nu+s\big),
\end{multline*}
where we used $c_{-s}=\frac{s}\ell$. Hence we derive $\sum_{i=1}^nx_i\nabla_i f_j=0$.
\end{proof}

\begin{remark}
In the case when the parameters satisfy $\ell | (s-q)$ the singular polynomials~\eqref{fi} appeared earlier in \cite[Proposition 4.1]{EtChm} where they were presented using the residue at infinity. In the other cases the space spanned by~\eqref{fi} does not contain the singular polynomials from \cite[Proposition 4.1]{EtChm} except the case $n=2$, $s=0$ and the case $n=1$.
\end{remark}

\vspace{5mm}

 {\bf Acknowledgements.} M.F. is very grateful to Y. Burman for stimulating discussions particularly in the beginning of the work, and to I.~Strachan for useful discussions and for pointing out the paper~\cite{Talamini}. M.F. is grateful to S.~Krivonos, O.~Lechtenfeld, A.~Marshakov and A.~Varchenko for stimulating discussions. Both authors are grateful to C.~Dunkl and S.~Griffeth for useful comments. The work of A.S. was partially supported by the RFBR grant 09-01-00239-a. The work of both authors was supported by  the EPSRC grant EP/F032889/1.

\end{document}